\documentclass[a4paper,12pt]{amsart}

\usepackage{amssymb,amsbsy,amsmath,amsfonts,amssymb,amscd}
\usepackage{latexsym}
\usepackage{amsxtra}
\usepackage{amscd}
\usepackage{graphics}
\usepackage{epic}
\usepackage{color}
\input xy
\xyoption{all}

\usepackage{ mathrsfs, amsfonts}

\newcommand\sD{{\mathcal D}}
\newcommand\sE{{\mathcal E}}

\newcommand\sF{{\mathcal F}}

\newcommand\sL{{\mathcal L}}

\newcommand\sN{{\mathcal N}}

\newcommand\sH{{\mathcal H}}
\newcommand\LL{{\mathbb L}}

\newcommand\om{\omega}

\newcommand\z{\zeta}

\newcommand\Ga{\Gamma}
\newcommand\De{\Delta}
\newcommand\ga{\gamma}
\newcommand\AL{\alpha}
\newcommand\bb{\beta}
\newcommand\e{\epsilon}
\newcommand\de{\delta}

\def\Bbb{\bf}

\newcommand{\CC}{\ensuremath{\mathbb{C}}}
\newcommand{\RR}{\ensuremath{\mathbb{R}}}
\newcommand{\ZZ}{\ensuremath{\mathbb{Z}}}
\newcommand{\QQ}{\ensuremath{\mathbb{Q}}}

\newcommand{\hol}{\ensuremath{\mathcal{O}}}

\newcommand{\HH}{\ensuremath{\mathbb{H}}}
\newcommand{\PP}{\ensuremath{\mathbb{P}}}

\newcommand{\HHH}{\ensuremath{\mathcal{H}}}

\newcommand{\ra}{\ensuremath{\rightarrow}}

\def\eea{\end{eqnarray*}}
\def\bea{\begin{eqnarray*}}

\newcommand{\Proof}{{\it Proof. }}

\newcommand\dual{\mathrel{\raise3pt\hbox{$\underline{\mathrm{\thinspace d
\thinspace}}$}}}
\newcommand\qe{\ifhmode\unskip\nobreak\fi\quad $\Box$}       

\def\BOX{\hfill\lower.5\baselineskip\hbox{$\Box$}}

\newcommand\R{\Bbb R}

\newtheorem{theo}[equation]{Theorem}
\newtheorem{remarkk}[equation]{Remark}
\newenvironment{rem}{\begin{remarkk}\rm}{\end{remarkk}}

\newtheorem{defin}[equation]{Definition}

\newtheorem{prop}[equation]{Proposition}
\newtheorem{cor}[equation]{Corollary}
\newtheorem{lemma}[equation]{Lemma}
\newtheorem{example}[equation]{Example}

\newtheorem{question}[equation]{Question}

\newcommand{\sR}{\ensuremath{\mathcal{R}}}

\newcommand{\Spec}{{\rm Spec\,}}

\newcommand{\To}{\;\longrightarrow\;}

\newcommand{\Mapsto}{\;\longmapsto\;}

\newcommand{\rk}{{\rm rk}}

\begin{document}

\title[Question on VHS by Fujita]{ Answer to a question by Fujita on Variation of Hodge Structures }
\author{ Fabrizio Catanese - Michael  Dettweiler}
\address {Lehrstuhl Mathematik VIII -Lehrstuhl Mathematik IV\\
Mathematisches Institut der Universit\"at Bayreuth\\
NW II,  Universit\"atsstr. 30\\
95447 Bayreuth}
\email{Fabrizio.Catanese@uni-bayreuth.de}
\email{Michael.Dettweiler@uni-bayreuth.de}

\thanks{
The present work took place in the realm of the DFG
Forschergruppe 790 ``Classification of algebraic
surfaces and compact complex manifolds'', the second named author  was supported by the DFG grant DE 1442/4-1}

\date{\today}

\maketitle
{\em  This article is  dedicated to  Yujiro Kawamata 
on the occasion of his $60$-th birthday.}

\begin{abstract}
We first provide details for the proof of Fujita's second  theorem for K\"ahler fibre spaces over a curve, asserting  that the direct image $V$ of the relative
dualizing sheaf splits as the direct sum $ V = A \oplus Q$, where $A$ is ample and $Q$ is unitary flat. Our main result then answers in the negative the question posed by Fujita
whether $V$ is semiample. In fact, $V$ is semiample if and only if $Q$ is associated to a representation of the fundamental group of $B$
having finite image. Our examples are based on hypergeometric integrals.

\end{abstract}

\tableofcontents
\section*{Introduction}

An important progress in classification theory was stimulated by  a theorem of Fujita, who showed (\cite{fuj1}) that if $X$ is a compact K\"ahler
manifold and $ f : X \ra B$ is a fibration onto a projective curve $B$ (i.e., $f$ has connected fibres), then the direct image sheaf
$$  V : = f_* \om_{X|B} = f_* ( \hol_X (K_X - f^* K_B))$$ is a semipositive (i.e., nef)  vector bundle on $B$, meaning that each quotient bundle $Q$
of $V$ has degree $\deg (Q) \geq 0$.

In the note \cite{fuj2} Fujita announced the following stronger result:

\begin{theo}{\bf (Fujita, \cite{fuj2})}\label{fuj2}

Let $f : X \ra B $ be  a fibration of a compact K\"ahler manifold $X$ over a projective curve $B$, and consider 
the direct image  sheaf $$ V : = f_* \om_{X|B} = f_* ( \hol_X (K_X - f^* K_B)).$$
Then $V$ splits as a direct sum  $ V  = A \oplus Q$, where $A$ is an ample vector bundle and $Q$ is a unitary  flat bundle.
\end{theo}

Fujita sketched the proof, but referred to a forthcoming article concerning the positivity  of the so-called local exponents.

After Fujita's articles, appeared then Kawamata's articles \cite{kaw0} \cite{kaw1}, which proved the conjecture $C_{n,1}$ (the subadditivity of Kodaira dimension for
 such fibrations, $ Kod (X) \geq Kod(B) + Kod (F)$, where $F$ is a general fibre) demonstrating the semipositivity also for 
the direct image of the higher powers of the relative dualizing sheaf 
$$W_m : =  f_* ( \om_{X|B}^{\otimes m}) = f_* ( \hol_X (m(K_X - f^* K_B))).$$
Kawamata's calculations are more directly related to Hodge theory, and especially a simple lemma, concerning
the degree of line bundles on  a curve whose metric grows at most  logarithmically around a finite number of singular points, played a crucial role for semipositivity.
Later Kawamata extended his result to the case where the dimension of the base variety of the fibration is greater than one(\cite{kaw2}).
\medskip

A first purpose of our article is to provide the missing details concerning the proof of the second theorem of Fujita, using Kawamata's lemma and some crucial estimates
given by Zucker (\cite{zucker}) for the growth of the norm of sections of the $L^2$-extension of Hodge bundles.

It is important to have in mind Fujita's second theorem in order to understand the question posed by Fujita in 1982 ( Problem 5, page 600 of \cite{katata},
Proceedings of the 1982 Taniguchi  Conference).
 
 \begin{question} {\bf (Fujita)}
 Is the direct image  $V : = f_* \om_{X|B}$ semi-ample ?
 
 \end{question}
 Saying that  a vector bundle $V$ is semi-ample means that the hyperplane divisor $H$ on $\PP : = Proj (V)$ is such that
 there exists a positive integer $m$ with $ | m H|$ base-point free.

 In our particular case, where $ V  = A \oplus Q$ with $A$ ample and $Q$ unitary  flat, it simply means that  the representation of the fundamental group 
 $ \rho : \pi_1(B) \ra  U (r, \CC)$ associated to the flat unitary
 rank-r bundle $Q$  has finite image (cf. theorem  \ref{semiample}).
 
 The main purpose of this article is to show that the question by Fujita has a negative answer.
 
 \begin{theo}\label{maintheorem}
There exist  surfaces $X$ of general type endowed with  a fibration $ f : X \ra B$ onto a curve $B$ of genus $\geq 3$, and with fibres of genus $6$,  such that $V : = f_* \om_{X|B}$ 
splits as a direct sum $ V = A  \oplus Q_1 \oplus Q_2$, where $A$ is an ample   rank-2  vector bundle, and the flat unitary rank-2 summands $Q_1, Q_2$
have infinite monodromy group (i.e., the image of $\rho_j$ is infinite). In particular, $V$ is not semi-ample.

 \end{theo}

 The surfaces in question are constructed using hypergeometric integrals associated to a cyclic group of order 7,
 and the non finiteness of the monodromy is a consequence of the classification
 due to  Schwarz (\cite{Schwarz}).
 
 A  minor contribution of the present paper is given also by the following result.
 
  \begin{theo}\label{Kodaira}
Let  $ f : X \ra B$  be a Kodaira fibration, i.e., $X$ is a surface and all the fibres of $f$ are smooth curves not all isomorphic to each other.
Then the direct image sheaf  $V : = f_* \om_{X|B}$ has strictly positive degree hence $\HHH : = R^1 f_* (\CC) \otimes \hol_B$ is a flat bundle which is not nef (i.e., not numerically semipositive).
 \end{theo}
 
  The underlying philosophy that theorems \ref{Kodaira} and \ref{semiample} convey is that the behaviour of flat bundles which are not unitary is rather wild.

\section{Preliminaries and reduction to the semistable case}

\subsection{Semipositive vector bundles on curves}

Let $B$ be a smooth complex projective curve, and assume that $V$ is a holomorphic vector bundle over it,
which we identify to its sheaf of holomorphic sections.

Recall the classical definition used by Fujita in \cite{fuj1}, \cite{fuj2}.

\begin{defin} 
Consider the projective bundle $\PP : = Proj (V) = \PP (V^{\vee})$, and the tautological divisor $H$ such that,
$ p : \PP \ra B$ being the natural projection, $p_* (\hol_{\PP} (H)) = V$.

Then $V$ is said to be:

(NP)  numerically semi-positive if and only if every quotient bundle $Q$ of $V$ has  degree $deg(Q) \geq 0$,

(NEF) nef if and only if $H$ is nef  on $\PP$,

(A) ample  if and only if $H$ is ample  on $\PP$,

(SA) semi-ample  if and only if a positive multiple of $H$ is spanned by global sections on $\PP$.

Recall that obviously ample implies semi-ample, semi-ample implies nef, while $H$ is nef if and only if $H$ is in the closure of the ample cone, or, 
equivalently, $ H \cdot C \geq 0$ for every effective curve $C$. 

As we shall recall, {\bf the conditions: nef and numerically semi-positive are equivalent}.

\end{defin}

\begin{rem}\label{quotient}
(1) Observe that if $U$ is a quotient bundle of $V$, then  $Proj (U)$ embeds in $ Proj (V)$ and the tautological divisor restricts to the tautological divisor,
hence if $V$ is ample (respectively, nef) then each quotient bundle $U$ is also ample (respectively, nef).
\end{rem}

We give an alternative proof\footnote{We became aware that a proof similar to ours is contained in Lazarsfeld's book \cite{laz}, theorem 6.4.15.} of a result of Hartshorne (Theorem 2.4, page 84 of \cite{hartVBC}), which is important  for our purposes.

\begin{prop}\label{ample}
A vector bundle $V$ on a curve   is  nef if and only it is numerically semi-positive, i.e.,  if and only if every quotient bundle $Q$ of $V$ has  degree $deg(Q) \geq 0$,
and $V$ is ample  if and only if every quotient bundle $Q$ of $V$ has  degree $deg(Q) > 0$.

\end{prop}

\Proof

One implication was essentially observed in greater generality in (1) of remark \ref{quotient}, except that we should show that a nef bundle has positive degree,
and an ample bundle has strictly positive degree.

By the Leray-Hirsch theorem the cohomology of $\PP$ is a free module over the cohomology of $B$, and its Chow ring is isomorphic to
$$ \ZZ[F,H] / ( H^r - d H^{r-1} F)$$

 where $r := rank (V) $, $ d : = \deg (V) = c_1 (V)$, and $F$ is a fibre of $ p : \PP \ra B$.
 
 By the same theorem, for every quotient bundle $Q$ of rank $k$ and degree $d'$ we obtain a projective subbundle $\PP' : = Proj (Q)$ such that
 the Chow ring of $\PP'$ equals $ \ZZ[F,H] / ( H^k - d' H^{r-1} F)$. 
 
 {\bf Step 1}: if $d' <0$, then $ H^k \cdot \PP' = d' < 0$, so $H$ is not  nef; similarly, if $d' = 0$, then $H$ is not  ample. 
 
  {\bf Step 2}: if $H$ is not nef, then there exists an irreducible  curve $ C' \subset \PP$ such that $H \cdot C' < 0$. The curve $C'$ cannot be contained in a fibre,
  since $H$ is ample on $F$, hence there exists a finite morphism from the normalization $C$ of $C'$ to $B$, $ f : C \ra B$.
  
  The pull-back of $V$, $ W : = f^* (V) $ has a quotient line bundle $L$ corresponding to the section $ C \subset W$, and $\deg (L) = H \cdot C' < 0$.
  
       {\bf Step 3}: Consider $f : C \ra B$ as in Step 2. Consider the Harder-Narasimhan filtration of $V^{\vee}$,
      $$ 0 \ra E_1 \ra E_2 \ra  \dots  E_h =  V^{\vee}.$$
      Here the slope of $E_1$ is  maximal for the associated graded bundle, i.e.,
      $ \mu (E_1) : = \deg (E_1) / \rk(E_1) \geq \mu (E_2/ E_1) \geq \mu (E_3/ E_2) \dots $,
      and all quotients $E_{j}/ E_{j-1}$ are stable.
      
      The pull back of a stable bundle is semistable, hence from step 2 we obtain an inclusion $L^{\vee} \ra f^* (V^{\vee}) $
      and therefore the slope of $L^{\vee}$, which is strictly positive if $H$ is not nef, is smaller or equal to the slope  $ \mu (E_1)$.
      Hence $V$ has a quotient bundle $  (E_1)^{\vee}$ with strictly negative degree.
  
    {\bf Step 4}:
    
    Let us work out the respective cones $\sE ff$  of effective curves, resp. $\sN ef$ of nef divisors for $\PP$. The latter  is a cone in the vector space $NS(\PP)$ with basis $H,F$, and it is the dual of the cone
    spanned by effective curves in the dual vector space $N_1(\PP)$ where we take as basis $ L : = F \cdot H^{r-2}$ (a line contained in a fibre)
    and $\Ga$, where $\Ga$ is a minimal section, i.e., such that $ \Ga \cdot H = : m $ is minimal (observe that $  m \geq 0$ if $H$ is net).
    
    We have thus: $$ L \cdot F = 0, L   \cdot H = 1,    \Ga \cdot  F = 1 ,  \Ga \cdot H =  m .$$ 
    
    The above formulae show that (since  the cone $\sE ff$ contains $L, \Ga$)  $F$, which is movable, hence nef, is nef but not ample;
    so $F$ is a boundary ray of $\sN ef$,
    while $L$ is a boundary ray for $\sE ff$.
    
    Assume that a curve $ a L + b \Ga$ is effective: then intersecting with $F$, which is movable, hence nef, we get $ b \geq 0$, and indeed $ b >0$ unless
    the curve is contained in a fibre (hence a multiple of $L$). Hence we may assume that the other boundary ray of $\sE ff$ is spanned by
    $ \Ga - a L$, where $ a \geq 0$. 
    
    Its orthogonal divisor class is given by $0 =  (x F + y H) \cdot  ( \Ga - a L) = x + y (m -a) \Leftrightarrow x = y (a-m)$, hence it is the class 
    is the class $ (a-m) F + H$.
    
    We get that $H$ is nef (respectively, ample) if and only if $$a-m \leq 0 \Leftrightarrow m \geq a$$ ( resp., $ m > a$).

     {\bf Step 5}: Assume now that   $H$ is nef and not ample: we want to conclude that  $\mu (E_1) =0$, hence concluding that there is a degree zero
     quotient $V \ra (E_1)^{\vee}$.
     
     By Step 4, we get that $a=m$ and that there are irreducible  curves $C'$  with class $ -\AL L + \bb \Ga, \bb > 0, \AL > 0$, as soon as $ m \bb > \AL$. 
     
     On the normalization $C$ of $C'$ we pull back via $ f : C \ra B$, and observe that $V$ has a line bundle quotient $L$ with slope
     $ C' \cdot H = - \AL + \bb m$. Again this slope, which is non negative, is bigger than the slope of $E_1$.  
     
     Take now the limit as $ \AL/  \bb  $ tends to $m$: then we conclude that  the slope of $E_1$ satisfies $\mu (E_1) \leq 0$; since $H$ is nef,
     we already know that $\mu (E_1) \geq 0$, hence $E_1$ has degree zero.

\qed

\begin{rem}
In general an extension $ 0 \ra W \ra V \ra E \ra 0$, where $W$ is ample, and $E$ is nef of degree zero, does not split.

Since the extension class lies in $ H^1 ( B, W \otimes E^{\vee})$, the dual space to $ H^0 ( B, E \otimes W^{\vee} \otimes K_B)$,
which is non zero if $B$ has genus $ g \geq 2$ and $\rk(E) =\rk(W) = \deg(W) = 1$.

\end{rem}

We give here a direct proof of the characterization of semi-ample unitary flat bundles; one step of the proof is related to a more general  theorem of Fujiwara (\cite{fujiwara}),
concerning semi-ample bundles with determinant of Kodaira dimension equal to zero.

\begin{theo}\label{semiample}
Let $\sH$ be a unitary flat vector bundle on a projective manifold $M$, associated to a representation $\rho : \pi_1(M) \ra U (r, \CC)$.
Then $\sH$ is nef and moreover $\sH$ is semi-ample if and only if $ Im (\rho )$ is finite.
\end{theo}

\Proof

Since $\sH$ is unitary flat,  $\sH$ is a Hermitian holomorphic bundle, and  by the principle   `curvature decreases in Hermitian subbundles' (page 79 of \cite{GH}, see also \cite{demailly} Prop. VII.6.10)  each subbundle has semi-negative degree and each quotient bundle $W$ of $\sH$ has semi-positive degree, hence $\sH$ is nef.

Assume that $\sH$ is semi-ample, and
let $B$ be a general linear curve section of $ M \subset \PP^N$, so that by Lefschetz' theorem we have a surjection $\pi_1(B) \ra  \pi_1(M)$.

Then $\sH | _B$ is also semi-ample and flat, corresponding to the composition homomorphism $ \pi_1(B) \ra  \pi_1(M) \ra U (r, \CC)$.

So, w.l.o.g., we may assume that $M$ is a curve $B$. 

{\bf Step I:}
we shall show that there exists a finite morphism $ p : B' \ra B$ such that the
pull back $ p^* (\sH)$ is a trivial holomorphic bundle. 

{\bf Step II:} a  unitary flat vector bundle on a projective curve is a trivial holomorphic bundle if and only if the associated representation
is trivial.

Step I and II , when put together with the following  lemma \ref{finite}, stating that the image of $\pi_1(B')$ has finite index 
in  $ \pi_1(B)$, suffice to show the difficult implication.  Since $\rho$ is trivial on the image of $\pi_1(B')$  by Step II, therefore the image  $ Im (\rho )$ is finite.

{\em Proof of Step I.}

Let $ \PP : = Proj (\sH)$ and let $ \pi : \PP \ra B$ the projection, and let  $F \cong \PP^{r-1}$ be a fibre. The hypothesis that $\sH$ is semi-ample means that there exists
a positive integer $ m \geq 1$ such that the linear system $ |mH|$ yields a morphism $ \psi : \PP \ra \PP^N$, which is finite on each fibre, since $\hol_F(H) = \hol_F(1)$.

We may choose $r$ divisors $ D_1, \dots, D_r \in | mH|$ such that $ D_1 \cap  \dots \cap  D_r \cap F = \emptyset.$ 
 
 Therefore we find $r$ multi-sections of $\pi$, setting $$ C_h : =   D_1 \cap  \dots \cap   \hat{D_h}  \dots  \cap D_r .$$
Let $B'$ be an irreducible component of the normalized fibre product  $ C_1 \times_B C_2 \times_B \dots \times_B C_r$:
then the pull back  $\sH '$ of $\sH$ admits $r$ sections of $ \hol_{\PP'}( H') $ yielding a birational map to $ B' \times \PP^r$.

Hence we get an injective homomorphism

$$  0 \ra \hol_{B'}^r \ra \sH ' \ra \sF \ra 0 $$
where the cokernel $\sF$ is concentrated on a finite set.

But then, since $ 0 = \deg (\sH ') = length (\sF)$, we obtain the desired isomorphism $ \hol_{B'}^r \cong \sH ' $.

{\em Proof of Step II.}

Let $B$ be a projective curve and  $\rho :  \pi_1(B) \ra   U (r, \CC)$ a unitary representation, and $\sH_{\rho}$ the associated 
flat holomorphic bundle. Since $\rho$ is unitary, it is a direct sum of irreducible unitary representations $\rho_j, j = 1, \dots k$.

Accordingly, we have a splitting 

$$ \sH_{\rho} = \oplus_{j=1}^k \sH_{\rho_j} .$$ 

Narasimhan and Seshadri have proven (see corollary 1, page 564 of \cite{NS}) that each $ \sH_{\rho_j}$ is a stable degree zero holomorphic bundle
on $B$.  Now, if $\rho$ is nontrivial, there exists an $h$ such that $\rho_h$ is non trivial. 

Assuming that we have an isomorphism $ \sH_{\rho} \cong \hol_B^r$ we derive a contradiction.

In fact, we have a surjection $$  \hol_B^r \ra \sH_{\rho_h} .$$

However, for each summand $  \hol_B$ its image must be trivial in $\sH_{\rho_h} $ if the latter bundle has rank $\geq 2$:
since its image would be a line bundle of the form $\hol_B(D)$, where $D$ is an effective or trivial divisor, and this would contradicit the stability
of $\sH_{\rho_h} $. But then the surjection $  \hol_B^r \ra \sH_{\rho_h} $ would be equalto zero, an obvious contradiction.

Assume instead that the rank of $\sH_{\rho_h} $ is equal to one. Then, since some summand would have a nonzero map $  \hol_B \ra \sH_{\rho_h} $,
this homomorphism would be an isomorphism since both line bundles have degree zero. Hence we would have $  \hol_B \cong \sH_{\rho_h} $,
contradicting theorem 2, page 560 of \cite{NS} (which asserts that  two  stable bundles corresponding to unitary representations
are isomorphic  if and only if the corresponding unitary representations are equivalent), since ${\rho_h} $ is nontrivial.

\medskip
{\em Proof of the easy implication.}

Conversely, if $ Im (\rho )$ is finite, there exists an \' etale Galois cover $ p : M' \ra M$, with Galois group $G$, such that $\sH' = p^* (\sH)$ is trivial.

We have $ \PP = (M' \times \PP^{r-1} )/ G$. For each point $x \in \PP^{r-1} $ we consider the $G$-orbit of $X$, and take a linear form  
$h'$ such that $h'$ does not vanish on  the orbit $ G x$: then the product of the $G$-transforms of $h'$ yields a section of $ \hol (m H')$
(here $ m = |G|$) 
which is $G$-invariant and does not vanish on $x$. Hence $\sH$ is semi-ample.

\qed

\begin{lemma}\label{finite}
Let $ p : B' \ra B$ be a finite morphism of curves. Then the image $\Ga$ of $\pi_1(B')$ has finite index 
in  $ \pi_1(B)$. 

\end{lemma}

\Proof
Let $B^*$ be the maximal open set such that, setting $ B'^* : = p^{-1} (B^*)$, $ p : B'^* \ra B^*$ is a finite unramified covering.

Then $\pi_1(B'^*)$ is a finite index subgroup of $\pi_1(B^*)$ and we conclude since $\pi_1(B^*)$ subjects onto $\pi_1(B)$,
and similarly $\pi_1(B'^*)$ subjects onto $\pi_1(B')$.

\qed 

\subsection{Semistable reduction}

Assume now that $ f : X \ra B$ is a fibration of a compact K\"ahler manifold $X$ over a projective curve $B$, and consider 
the invertible sheaf $ \om : = \om_{X|B} = \hol_X (K_X - f^* K_B).$ 

By Hironaka's theorem there is a sequence of blow ups with smooth centres $ \pi: \hat{X} \ra X$ such that $$ \hat{f} := f \circ \pi: \hat{X} \ra B$$
has the property that all singular  fibres $F$ are such that $ F = \sum_i m_i F_i$, and $F_{red} =  \sum_i  F_i$ is a normal crossing divisor.

Since $\pi_* \hol_{\hat{X} } (K_{\hat{X} }) = \hol_X (K_X)$ we obtain 

$$ \hat{f} _* \om_{\hat{X}|B } = \hat{f} _* \hol_{\hat{X} } (K_{\hat{X} } - \hat{f} ^* K_B ) = f_* \hol_X (K_X - f^* K_B) = f_*  \om_{X|B} .$$

Therefore we shall assume wlog that all the reduced fibres of $f$ are normal crossing divisors.

\begin{theo} {\bf (Semistable reduction theorem, \cite{toroidalEmb})}
There exists a cyclic Galois covering of $B$, $B' \ra B = B' /G$, such that the normalization $X''$ of the fibre product $ B' \times_B X$
admits a resolution $ X' \ra X''$ such that the resulting fibration $f' : X' \ra B'$ has all the fibres which are reduced and 
normal crossing divisors.
\end{theo}

\begin{equation*}
\xymatrix{
X' \ar[d]^{f'}\ar[r]^{v'} & X'' \ar[d]\ar[r]^{v''}  & X\ar[d]^f\\
B' \ar[r]^{Id} &B'  \ar[r]^u & B ,}
\end{equation*}

\begin{rem}
At each singular fibre $ F = \sum_i m_i C_i$ corresponding to a point $t=0$ on $B$, the theorem yields   a base change $ t = \tau^n$,
where $ m_i | n, \forall i$, and $ n >> 0$.

As a notation, we set $v := v'' \circ v' : X' \ra X$. We set also $ n = m_i d_i$.

\end{rem} 

\begin{prop}\label{sstablered}
The sheaf $V' : =   f'_* \om_{X'|B'}$ is a subsheaf of the sheaf $ u^* (V)$, where $V : =   f_* \om_{X|B} $, and the cokernel 
$ u^* (V) / V'$ is concentrated on the set of points corresponding to singular fibres of $f'$.

In particular, since  $V$ and $V'$ are semipositive by Fujita's first theorem, if $V'$ satisfies the property that for each degree  0 
quotient bundle $Q'$ of $V'$ then there is a splitting
$ V' = E' \oplus Q' $ for the projection $ p : V' \ra Q'$,and $Q'$ is unitary flat,  then $V'$ splits as the direct sum
$ V'  = A \oplus Q$, where $A$ is an ample vector bundle and $Q$ is a  flat unitary bundle, and the same conclusion holds also for $V$.

\end{prop}

\Proof
It suffices to work locally around each point $ P' \in B'$, mapping to a point $P$ corresponding
to a singular fibre $F$ of $f$ and consider the base change $ t = \tau^n$, where $n $ 
may be assumed not to depend on the point $P$. 

By the Hurwitz formula $$K_{B'} = u^* (K_B) + \sum_{P'} (n-1)P' , \ \   K_{X'} = v^* K_X + R,$$
hence
$$ \om_{X'|B'}=  K_{X'}  - f'^* (K_{B'} ) = v^* (\om_{X|B}) -  ( \sum_{P'}  (n-1) F' - R),$$
and our assertion would be  proven if the divisor $ \sum_{P'}  (n-1) F' - R$, supported on the inverse images of the singular fibres of $f$, 
is effective. 

Let us  work locally around the fibre $F'$ of $f'$ which lies above the point $P'$.

Recall the following lemma, where (see \cite{km}, Lemma 5.12, page 156), equality holds, for $X''$ Cohen Macaulay,
and $X'$ the resolution, if $X''$ has rational singularities.

\begin{lemma}
Let $g : X' \ra X''$ be a birational morphism between normal varieties: then $g_* (\hol(K_{X'} )) \subset  \hol(K_{X''} )$.
\end{lemma}

The first thing to do is to separate in $R$ the $v$-exceptional divisors and the divisors $D_i$, which are the strict transforms of $F_i$.
Recall in fact that, if $\ga_i = 0$ is a local equation of $F_i$, then
$$\tau^n =  \ga_i ^{m_i} \Leftrightarrow \Pi_ {\epsilon^{m_i} = 1}( \tau^{d_i} - \epsilon \ga_i),$$
and the local equation of $D_i$ in the normalization $X''$ becomes $\tau = 0$.

Therefore $ R \geq \sum_i (d_i -1) D_i$.

Finally, we get, working again locally, and observing that $ n - d_i = d_i (m_i -1) \geq 0$:

$$ u_*   f'_* \om_{X'|B'} \subset  u_*   f''_* \om_{X''|B'} =  u_*  f''_*  (v''^* \om_{X|B} + \sum_i (d_i -1 - (n-1)) D_i ) = $$ $$=u_*   f''_*  (v''^* \om_{X|B} -  \sum_i (m_i d_i - d_i ) D_i )\subset f_* \om_{X|B}. $$

We are left to prove the second assertion of proposition \ref{sstablered}. For this purpose we consider again the Harder-Narasimhan filtration of $V$,
      $$ 0 \ra V_1 \ra V_2 \ra  \dots  V_h =  V,$$
      where as usual the slope of $V_1$ is maximal. It gives rise to an exact sequence $  0 \ra  U  \ra V \ra   Z \ra 0$,
      where $U$ is ample and $ \deg(Z) = 0$.
      
      We have then a generically invertible homomorphism between two vector bundles of the same rank:
      $$ V '  \ra  u^* (V) \supset  u^*(U )   \supset  0. $$

We set $ Q : = u^* (Z)$, and observe that  $V'  \ra Q = u^* (Z)$ must be  surjective, else $V'$ would have a negative degree quotient. Then, by our assumption, it follows that we
have a splitting $ V' = E' \oplus Q$.

We claim that the inclusion $ Q \subset V'$ induces a splitting of  $  0 \ra  U  \ra V \ra   Z \ra 0$, yielding
$ V = U \oplus Z$. This follows since $ V = u_*(u^*V)^G$, hence the homomorphism $ V' \ra u^*V$ induces a chain of 
homomorphisms $$ Z = u_*(u^*Z)^G  \ra  u_*(u^*Z)^G \oplus   u_*(A)^G  = u_*(V')^G \ra V \ra Z $$
whose composition is the identity.

We show that $Z$ is a flat bundle.

Since $Q$ is a flat unitary bundle, $Q$ is a quotient $( \tilde{B'} \times \CC^r)/\pi_1(B')$, where $\tilde{B'}$ is the universal covering of $B'$.
The action is determined by a homomorphism $ \rho' : \pi_1 (B') \ra  U (r, \CC)$.

Denote by $B^*$ the complement of the branch locus of $ B' \ra B$, and by $B'^*$ its inverse image.
Since $ Z = u_* (Q)^G$, we have that the restriction $ Z | _{B^*}$ is a flat bundle associated to a representation $ \rho^* : \pi_1 (B^*) \ra  U (r, \CC)$.

However, since $Z$ is a vector bundle on $B$, the restriction of $ \rho^* $ to generators of the kernel of the surjection $ \pi_1 (B^*) \ra \pi_1 (B) $
is trivial, hence $ \rho^* $  factors through  $ \rho : \pi_1 (B) \ra  U (r, \CC)$.

Since the restriction of $Q$ to $B'^*$ corresponds to the restriction of $ \rho^* $ to  $ \pi_1 (B'^*)$, and it is trivial on the kernel 
of $ \pi_1 (B'^*) \ra \pi_1 (B') $, we have shown that $\rho'$ factors through $\rho$.

It follows that $Z$ is a flat bundle: we have in fact seen that 
 $Q$ is a quotient $( \tilde{B'} \times \CC^r)/\pi_1(B')$, where $\tilde{B'}$ is the universal covering of $B'$, and where 
the action is determined by the homomorphism $ \rho' : \pi_1 (B') \ra U(r, \CC)$.

Hence $Z$ is a quotient 
$( \tilde{B'} \times \CC^r)/ \pi_1^{orb}(B' \ra B)$, where the orbifold fundamental group is defined (see e.g. \cite{cime}, pages 101 and following for more details)
by the extension
$$ 1 \ra \pi_1(B') \ra  \pi_1^{orb}(B' \ra B) \ra G \ra 1 ,$$ such that $B$ is the quotient of $\tilde{B'}$ by $ \pi_1^{orb}(B' \ra B)$.

$ \pi_1^{orb}(B' \ra B)$ is a quotient of the fundamental group of  $B^*$, the complement of the branch locus of $ B' \ra B$. 

We also saw that since   $Q$ is the  pull back of $Z$,
the representation on $\CC^r$ of the orbifold fundamental group $\pi_1^{orb}(B' \ra B)$ factors through the surjection
$\pi_1^{orb}(B' \ra B) \ra \pi_1 (B)$, therefore  $Z = ( \tilde{B}  \times \CC^r)/\pi_1(B)$ is a flat unitary bundle over $B$.

\qed

\section{Fujita's second  theorem}

In this section we shall use some standard differential geometric terminology which we now recall.

\begin{defin}
Let $(E,h)$ be an Hermitian vector bundle on a complex manifold $M$. Take the canonical Chern connection 
associated to the Hermitian metric $h$, and denote by $\Theta(E,h)$ the associated Hermitian curvature,
which gives a Hermitian form on the complex vector bundle bundle $T_M \otimes E$.

Then (see for instance  \cite{laz}, and also \cite{kob}), one says that $E$ is Nakano positive (resp.: semi-positive) if there exists a Hermitian metric $h$
such  that the Hermitian form associated to $\Theta(E,h)$ is strictly positive definite (resp.: semi-positive definite).

In local coordinates $(z_1, \dots, z_m)$ there exists a local frame $e_{\lambda}$ such that

$$ i   \Theta(E,h) = \Sigma_{j,k,\lambda, \mu} c_{j,k,\lambda, \mu} dz_j \wedge d \overline{z_k} \otimes e_{\lambda}^* \otimes e_{\mu}, 
\ c_{k,j, \mu, \lambda} =  \overline{c_{j,k,\lambda, \mu} }.$$  

While one says that $E$ is Griffiths positive (resp.: semi-positive) if there exists a Hermitian form as above which  is positive on rank 1 tensors 
$T_M \otimes E$. Nakano positive implies Griffiths positive, Griffiths positive implies ample, and Griffiths semipositivity implies nefness.

If $M$ is a curve, then 
$$ i   \Theta(E,h) = \Sigma_{\lambda, \mu} c_{\lambda, \mu}   e_{\lambda}^* \otimes e_{\mu}  \otimes  dz \wedge d \overline{z},$$ 
and Nakano and Griffiths positivity (resp. : semi-positivity) coincide, since they both boil down  to the requirement that the Hermitian matrix 
$(c_{\lambda, \mu})$ is positive definite (resp. : semi-positive), and we shall simply then say that an Hermitian vector bundle is positive
(resp. : semi-positive). These notions then imply respectively ampleness and numerical semi-positivity (nefness) of the bundle $E$.

\end{defin}

\begin{rem}
Umemura proved (\cite{um}, theorem 2.6, see also \cite{CF})  that a vector bundle $V$ over a curve $B$ is positive (i.e.,  Griffiths positive, or equivalently Nakano positive) if and only if $V$ is ample.
\end{rem}

We pass now to Fujita's second theorem.

\begin{theo}{\bf (Fujita, \cite{fuj2})}\label{fuj2}

Let $f : X \ra B $ be  a fibration of a compact K\"ahler manifold $X$ over a projective curve $B$, and consider 
the direct image  sheaf $$ V : = f_* \om_{X|B} = f_* ( \hol_X (K_X - f^* K_B)).$$
Then $V$ splits as a direct sum  $ V  = A \oplus Q$, where $A$ is an ample vector bundle and $Q$ is a unitary flat bundle

\end{theo}

The details of the proof were never published by Fujita. Thanks to the auxiliary results shown in the previous section,
in particular proposition \ref{sstablered},
it suffices to prove the theorem in the semistable case, i.e., where each fibre is reduced and a normal crossing divisor.

\Proof
We first treat  the case where there are no singular fibres, and the underlying idea is simpler.

{\bf Case 1: there are no singular fibres}

In this case $V$ is semipositive, as it was shown by Fujita in \cite{fuj1}.

Another proof via Hodge bundles was given by Griffiths in \cite{griff2} (see also \cite{Griffiths_Topics} and \cite{zucker-remarks}).

The underlying idea runs as follows.

$V$ is a holomorphic subbundle of the holomorphic  vector bundle $\sH$  associated to the local
system $\HH : = \sR^n f_* (\ZZ_X)$, i.e.,  $\sH = \HH \otimes_{\ZZ} \hol_B$.

In fact, we have that  $ \bar{V}$ is an antiholomorphic subbundle and $V \oplus \bar{V}  \subset \sH$ is a subbundle such that 
 the  Hermitian orthogonal splitting  $V \oplus \bar{V}$  identifies
$ \bar{V}  $ to the dual bundle $V^{\vee}$. The bundle $\sH$ is flat, hence the curvature $\Theta_{\sH} $ associated to the flat connection satisfies 
$\Theta_{\sH} \equiv 0$.
(in particular, see \cite{kob}, proposition 3.1 (a), page 42, all the real Chern classes of a flat bundle vanish).

We  view $V$ as a holomorphic subbundle of $\sH$, while  $$V^{\vee} \cong R^nf_* \hol_X, \ n = \dim(X) - 1$$ is a holomorphic quotient bundle of $\sH$.
\footnote{It is important to remark that we take here the curvature of a flat, but not unitarily flat, bundle $\sH$; in particular the principle:
curvature decreases in subbundles (page 79 of \cite{GH}) does not hold, since this assumes that we take the curvature
associated to an Hermitian metric, while the intersection form on $\HH$ is not definite.}

Using arguments similar to  the curvature formula for subbundles (see \cite{Griffiths_Topics}, Lecture 2)
$$  \Theta_{V} =  \Theta_{\sH} |_V + \bar{\sigma} \ ^t \sigma = \bar{\sigma} \ ^t \sigma ,$$
Griffihts proves (\cite{griff2}, see also corollary 5, page 34 \cite{Griffiths_Topics}) that the curvature of $V^{\vee}$ is semi-negative, 
since its local expression is of the form
$ i h' (z) d \bar{z} \wedge dz $, where $h'(z)$ is a semi-positive definite Hermitian matrix. 
In particular we have that  the curvature  $  \Theta_{V} $  of $V$
is semipositive and, 
moreover,  that the curvature  vanishes identically  if and only if the second fundamental form $\sigma$
vanishes identically, i.e., if and only if $V$ is a flat subbundle.

However, by semi-positivity, we get that  the curvature  vanishes identically  if and only its integral, the degree of $V$,
equals zero. Hence $V$ is a flat bundle if and only if  it has degree $0$.

The same result then  holds true, by an identical reasoning,  for each holomorphic quotient bundle $Q$ of $V$ by the following argument.

Assume now that $V$ is not ample. By Hartshorne's theorem (proposition \ref{ample}), there is an exact sequence of holomorphic bundles

$$ 0 \ra E \ra V \ra Q \ra 0 $$
where the quotient bundle $Q$ has degree $0$.

Dualizing, we obtain 

$$ 0 \ra Q^{\vee}  \ra V^{\vee}   \ra E^{\vee}  \ra 0 $$

and, since $V^{\vee}$ has semi-negative Hermitian curvature,
then by the cited principle `curvature decreases in Hermitian subbundles' (page 79 of \cite{GH}, see also \cite{demailly} Prop. VII.6.10) 
 $Q^{\vee}$ has semi-negative Hermitian curvature.

However, $Q^{\vee}$ has degree $0$, thus  the integral of the semi-negative curvature of $Q^{\vee}$  is zero,
so its Hermitian  curvature $\Theta_{ Q^{\vee}} \equiv 0$, hence $Q^{\vee} \cong \bar{Q} \subset \bar{V}$   is a flat subbundle of the flat bundle $\sH$,
 and similarly  $Q$ is a flat bundle.

Since we have an inclusion $ \bar{Q} \subset \bar{V} \subset \sH$ of the    flat antiholomorphic subbundle $ \bar{V}$, we obtain 
by complex conjugation  an inclusion of the holomorphic
subbundle $ Q \subset V $ hence  a splitting of the surjection $ V \ra Q$. 

Finally, $Q$ is unitary flat since the intersection form on $V$ is, up to constant, strictly positive definite.

\qed

{\bf Case 2: there are  singular fibres, which are normal crossing divisors, and the local monodromy is unipotent,
since the fibres are reduced.}

The treatment of the general case is similar: it suffices to show that the degree of $Q$ is the integral
of the curvature form on $ B^*$, where  $ B \setminus B^*= : S$ is the set of critical values of $f$.

Recall that the degree of the bundle $Q$ is the degree of its top exterior power, the so-called
determinant bundle $ \det(Q)$.

We use here a well known lemma (see lemma 5, page 61 of \cite{kaw1}, also proposition 3.4, page 11 of  \cite{Peters}):

\begin{lemma}
Let $L$ be a holomorphic line bundle over a projective curve $B$, and assume that 
$L$ admits a singular metric $h$ which is regular outside of a finite set $S$ and  has at most logarithmic growth at
the points $ p \in S$ (i.e., if $z$ is a local coordinate at $p$, then $ | h(z)| \leq C  log |z|^{-m }$, where $C$ is a positive constant,
and $m$ is a positive integer).

Then the first Chern form $c_1(L,h) := \Theta_h $ is integrable on $B$, and its integral equals $ deg(L)$.
\end{lemma}

Let us briefly recall how the existence of such a metric is shown to exist.

We have the VHS (variation of Hodge structure)  on the punctured curve $B^*$ given by the local system $$\HH^*   : = \sR^n F_* (\ZZ_{X^*}),$$
where $X^* : = f^{-1} (B^*)$ and $ F : X^* \ra B^*$ is the restriction of $f$ to $X^*$.

Again $V^* : = V | _{B^*}$ is a subbundle of the flat bundle $\sH^* : = \HH^* \otimes_{\ZZ} \hol_{B^*}$,
and we get a subbundle  $ V^* \oplus \bar{V^* } \subset \sH^* .$ 

$\sH^*$ is a flat holomorphic bundle and the associated holomorphic connection
$\nabla^*$ on $\sH^*$ is the so called Gauss-Manin connection.

We then have the Deligne canonical extension $(\sD \sH , \nabla)$ of the pair $(\sH^*, \nabla^*)$ to a holomorphic vector bundle 
$ \sD \sH$ endowed with a meromorphic connection $\nabla$ having simple poles on the points of $S$, and with nilpotent residue matrices.
We refer to part II of \cite{kollar} (see especially section 2 , and theorem 2.6) for more details about the  presentation of this extension, 
which we now briefly describe.

We let $D$ be the normal crossing divisor $ f^{-1} (S)$, and consider the relative De Rham complex 
$$ \Omega^{\cdot} _{X|B} (\log  D)$$
with logarithmic singularities along $D$.  

The hypercohomology sheaf $\sD \sH  =  \R^i f_* ( \Omega^{\cdot} _{X|B} (\log  D))$   gives an extension of $ \sH_i^* = \R^i F_* \CC_{X^*} \otimes \hol_{B^*}$
from $B^*$ to $B$, the Deligne extension.

By the work of Schmid (\cite{schmid}, see also \cite{G-S})  the Hodge filtration on $\sH^*$ extends to a holomorphic (decreasing) filtration $\sF^i (\sD \sH), \ i=0, \dots, 
n := dim X - 1.$

In particular, for $ i = n := dim X - 1$ we have, as proven by Kawamata in \cite{kaw1} (lemma 1, page 59)
$$ V = f_* \om_{X|B}  =   \sF^n (\R^n f_* ( \Omega^{\cdot} _{X|B} (\log  D))). $$

As explained in \cite{G-S}, and also in proposition 4.4., page 433 of \cite{zucker},  logarithmic forms are precisely those holomorphic forms
with  the property of being square integrable, 
and this approach was taken up by Zucker \cite{zucker} who used the explicit description of the limiting Hodge structure
found by Schmid in \cite{schmid} in order to prove the following result (which is proven in the course of the proof of proposition 5.2, pages 
435-436).

\begin{lemma}
For each point $s \in S$ there exists a  basis of $V$ given by elements $\sigma_j$ such that their norm 
in the flat metric outside the punctures grows at most logarithmically.

In particular, for each quotient bundle $Q$ of $V$ its determinant admits a metric with growth at most logarithmic at the punctures $ s \in S$,
and the degree of $Q$ is given by the integral of the first Chern form of the singular metric.

\end{lemma}
\begin{rem}
Observe that a similar result, but for the determinant of $V$, is used in \cite{Peters} to show the flatness of $V$ in the case that $ \deg (V) = 0$.
\end{rem}

Therefore we can conclude that, since $ \deg (Q) = 0$, and since its integral is given by the norm of the second fundamental form,
which is semipositive, then the  the second fundamental form vanishes identically and $Q$ is a flat sub-bundle.  The same  argument as 
 the one given for case 1 shows then that we have an inclusion $ Q^* \ra V^* : = V |_{B^*}$. 

Now  $Q^* : = Q _{| B^*}$ is a unitary  flat subbundle of the flat bundle $\sH^*$, in particular the local  monodromies at the punctures (the points of $S$),
being unitary and unipotent, are trivial:
hence   $Q^*$ has a flat extension to $B$ which we denote by $\hat{Q}$.

Clearly we have inclusions 
$$  \hat{Q} \subset V \subset \sD \sH,$$
and we obtain a homomorphism  
 $ \psi : \hat{Q}\ra Q$  composing the inclusion $  \hat{Q} \ra V$ with the surjection $ V \ra Q$.
 
From the fact that  $\psi$ is an isomorphism 
over $B^*$ we infer that $\psi$ is an isomorphism: since $\det (\psi)$ is  not identically zero, and is a section of a degree zero line bundle.

Hence we conclude that the composition of $\psi^{-1}$ with the   inclusion $  \hat{Q} \ra V$  gives then the desired splitting of the 
surjection $ V \ra Q$.

\qed

\begin{cor}\label{cor}

Let $f : X \ra B $ be  a fibration of a compact K\"ahler manifold $X$ over a projective curve $B$, and consider 
the direct image  sheaf $$ V : = f_* \om_{X|B} = f_* ( \hol_X (K_X - f^* K_B)).$$
Then $V$ splits as a direct sum  $ V  = A \bigoplus (\oplus_{i=1}^h Q_i)$, where $A$ is an ample vector bundle and each $Q_i$ is flat vector bundle
without any nontrivial degree zero quotient. 
Moreover,

(I) if $Q_i$ has rank equal to 1, then it is a torsion bundle ($\exists \  m$ such that $Q_i^{\otimes m}$
is trivial),

(II) if the genus of the  curve $B$ equals 1, then each $Q_i$ has rank 1.

(III)  In particular,  if the genus of the  curve $B$ is at most  1, then $V$ is semi-ample.

\end{cor}

\Proof
Each time $V$ has a degree zero quotient, this yields a splitting, as shown in theorem \ref{fuj2}.
Therefore, we obtain that  $V$ splits as a direct sum  $ V  = A \bigoplus (\oplus_{i=1}^h Q_i)$,
where $A$ is an ample vector bundle and each $Q_i$ is flat vector bundle
without any nontrivial degree zero quotient. 

(I) This was proven by Deligne in \cite{TH2}, cor. 4.2.8 (iii) (b).

(II) This is immediate, since the fundamental group of a curve $B$ of genus $1$ is abelian,
hence every representation splits as a direct sum of 1-dimensional representations.

(III) A torsion line bundle is semi-ample, and a direct sum of semi-ample vector bundles is semi-ample.

\qed

\begin{rem}
Part (III) of the above corollary was proven by Barja in \cite{Barja}.

\end{rem} 

We proceed now to prove 

{\bf Theorem \ref{Kodaira}}
{\em Let  $ f : X \ra B$  be a Kodaira fibration, i.e., $X$ is a surface and all the fibres of $f$ are smooth curves of genus
$g \geq 2$ not all isomorphic to each other.
Then the direct image sheaf  $V : = f_* \om_{X|B}$ has strictly positive degree hence $\HHH : = R^1 f_* (\CC) \otimes \hol_B$ is a flat bundle which is not nef (i.e., not numerically semipositive).}

\Proof
Since all the fibres of $f$ are smooth, we have an exact sequence 
$$  0 \ra V \ra \HHH \ra V^{\vee} \ra 0 ,$$
and it suffices to show that the degree of the quotient bundle $V$  is strictly negative, or, equivalently,
$ \deg (V) > 0$.

We have that $$ \deg (V) = K_X^2 - 8 (g-1)(g-1),$$ where $g$ is the genus of the fibres of $f$, and $b$ is the genus of $B$.
As well known (see \cite{bpv}) also the genus $ b \geq 2$, and therefore $X$ contains no rational curve and is therefore
a minimal surface.

Since $f$ is a differentiable fibre bundle, we have for  the Euler- Poincar\'e characteristic  of $X$
$$ e (X) = 4 (b-1) (g-1).$$

Kodaira (\cite{kod}) proved that for such fibrations the topological index $\sigma(X)$, the signature of the intersection form
on $H^2(X, \RR)$ is positive. By the index theorem (see again \cite{bpv}) we have
$$0 < 3  \sigma(X) = c_1^2(X) - 2 c_2(X) =  K_X^2 - 2 e(X) = K_X^2 - 8 (g-1)(g-1) = \deg (V).$$ 

\qed

\begin{cor}\label{exist}
There are flat bundles on curves which are not nef, in particular do not admit an Hermitian metric
with  semipositive curvature.

\end{cor}

\Proof
See b) of the following remark.

\qed

\begin{rem}
a) The examples of Kodaira (\cite{kod} and other examples (\cite{cat-rollenske}) show that the direct image $V$
does not need to be ample, since  it can have a trivial summand.

b) These and other examples of course,  showing the existence of Kodaira fibrations,  furnish the proof of 
corollary \ref{exist}.

\end{rem}

\section{A curve of genus 6 with cyclic symmetry of order 7}

In this section we explain how we obtain explicit examples of fibrations where $V = f_* \omega$ has a flat summand.

Consider the equation 
$$ z_1^7 = y_1 y_0 (y_1 - y_0 )  (y_1 - x y_0 )^4 , \ \  x \in \CC \setminus \{0,1  \}$$
describing a singular curve inside the weighted projective space $\PP (1,1,7)$,
with variables $ y_0,  y_1, z_1$ (alternatively, a curve inside the line bundle $\LL$  
over $\PP^1$ whose sheaf of holomorphic sections $\sL_1$ equals $ \hol_{\PP^1}(1)$.

Denote by $C$ the normalization of the above curve. $C$   has a Galois cover  $ \phi : C \ra \PP^1$ with Galois group 
$ G = \mu_7 = \{ \e | \e^7 = 1\}$, acting by $ z_1 \mapsto \e z_1$; there are exactly 4 ramification  points,
$P_0$ lying above $ y_0 = 0$, $P_1$ lying above $ y_1 = 0$, $P_2$ lying above $ y_1 - y_0 = 0$,
$P_x$ lying above $ y_1 - x y_0 = 0$.
Correspondingly there are nonzero sections $ w_0 \in H^0 (\hol_C (P_0))$,$ w_1 \in H^0 (\hol_C (P_1))$,
$ w_2 \in H^0 (\hol_C (P_2))$, $ w_x \in H^0 (\hol_C (P_x))$ such that we obtain a factorization
$$z_1 = w_0 w_1 w_2 w_x^4.$$
We shall alternatively use the classical notation $ w_0 = y_0^{\frac{1}{7}}$,$ w_1 = y_1^{\frac{1}{7}}$,
$ w_2 = (y_1 - y_0)^{\frac{1}{7}}$,$ w_x = (y_1 - x y_0)^{\frac{1}{7}}$.

Since the ramification points are $G$-invariant, the corresponding sheaves admit a $G$-linearization, and we shall
choose the linearization by which  the generator $g \in G$ acts via
$$  g(z_1)  =  \z z_1 , g (w_0) = w_0, g (w_1) = w_1, g (w_2) = w_2, \ g (w_x) =  \z^2 w_x,  $$
where $\z : = \exp (\frac{2 \pi i }{7} )$.

The biregular structure of $C$ is described (see \cite{par}, \cite{bc}, \cite{cyclic}) as
$$ C = \Spec (  \hol_{\PP^1} \bigoplus (\oplus_{j=1}^6 z_j \sL_j^{-1}).$$

Here, $ z_j \in H^0 (C, \phi^* (\sL_j))$ and $z_j \sL_j^{-1}$ is the $j$-th character eigensheaf, i.e., $ g (z_j) = \z^j z_j$.

 It is easy to describe the sections  $z_j$ by taking powers of the above equation for $z_1$ and reducing the exponents modulo $7$:
 $$z_2 = w_0^2 w_1^2 w_2^2 w_x, \ z_2 = w_0^3 w_1^3 w_2^3 w_x^5, \ z_3 = w_0^4 w_1^4 w_2^4 w_x^2,$$
 and then we observe that $$ z_j \cdot z_{7-j} = \de : =  y_1 y_0 (y_1 - y_0 )  (y_1 - x y_0 ).$$
 
 Hence we derive 
 $$\sL_1 = \sL_2 =   \hol_{\PP^1}(1) , \sL_3 =  \sL_4 =  \hol_{\PP^1}(2) , \sL_5 =  \sL_6 =  \hol_{\PP^1}(3).$$

The formulae by Pardini, proposition 4.1 (page 207 of \cite{par}) yield, when we denote by $\De$ the reduced branch divisor $ \de = 0$: 
$$( \phi_* \Omega^1_C)_j \cong  \Omega^1_{\PP^1} (\De) \otimes \sL_j^{-1} \cong  \Omega^1_{\PP^1} \otimes \sL_{-j},$$
and we have more precisely, using affine coordinates where $y_0 =1 , y : = y_1$ and setting $ \eta_\De : = \frac{dy}{\de}$ 

$$ V_j : = H^0 ( \Omega^1_C)_j = z_j \eta_\De \phi^* H^0( \PP^1, \sL_{-j} (-2)). $$

Hence $ V_5 = V_6 = 0$, while $ \dim V_3 = \dim V_4 = 1$, and finally
$$  V_1 = \{ (a_0 + a_1 y) z_1 \eta_\De | a_0, a_1 \in \CC \} =  \{ (a_0 + a_1 y) w_1^{-6}  w_2^{-6} w_x ^{-3}dy   \} $$
$$  V_2 =  \{ (b_0 + b_1 y) w_1^{-5}  w_2^{-5} w_x ^{-6} dy   | b_0, b_1 \in \CC\} .$$

In other words, a basis of $V_1$ is given by $\tau, y \tau$, where
$$ \tau : = y^{- \frac{6}{7}}  (y-1)^ {- \frac{6}{7}}  (y- x) ^{- \frac{3}{7}}  dy .$$

\begin{rem}
(I) Observe that changing the generator of $G$ with its opposite has the effect of replacing $\tau$ with 
$$ \tau' : = y^{- \frac{1}{7}}  (y-1)^ {- \frac{1}{7}}  (y- x) ^{- \frac{4}{7}}  dy .$$

(II) The curve $C$ has genus $6$, and the linear subsystems of the canonical system corresponding to the eigensheaves
have a base locus, since the greatest common divisor of the  elements in $V_1$ is  $w_x^3$, for $V_2$ it is $w_0 w_1 w_2$,
while $ V_3 = \{ (w_0 w_1 w_2)^2 w_x^4 \}$, $ V_4 = \{ (w_0 w_1 w_2)^3 w_x \}$.

\end{rem}

Consider now the Hodge decomposition of the cohomology of $C$, viewed as a $G$-representation:

$$ \HH : = H^1 (C, \CC) =  H^0 ( \Omega^1_C) \oplus \overline {H^0 ( \Omega^1_C)}, $$ 
$$ \HH =  (V_1 \oplus V_2  \oplus V_3 \oplus V_4) \oplus \overline { (V_1 \oplus V_2  \oplus V_3 \oplus V_4) } .$$

The consequence is that $$ \HH_j  =  V_j ,  \HH_{7-j}   = \overline { V_j } \ j = 1,2 . $$
while  $$ \HH_3  =  V_3 \oplus \overline { V_4 },$$
and similarly for $ \HH_4$.

We conclude the above discussion with its consequence
\begin{prop}\label{7}
Let $f : X \ra B$ be a semistable  fibration of a surface $X$ onto a projective curve, such that the group $G = \mu_7$
acts on this fibration inducing the identity on $B$. Assume that the general fibre $F$ has genus 6 and that
$G$ has exactly 4 fixed points on $F$, with tangential characters $ (1,1,1,4)$.

Then if we split $ V = f_* (\omega_{X|B})$ into  eigensheaves, then the eigensheaves 
$V_1, V_2$ are unitary flat rank 2 bundles.
\end{prop}
\Proof
Since the fibration is semistable, the local monodromies are unipotent: on the other hand, they are unitary,
hence they must be trivial. This implies that the local systems $\HH_1^*$ and $\HH_2^*$ 
 have  respective flat extensions to local systems  $\HH_1$ and $\HH_2$
on the whole curve $B$. Denote by $\sH_j : = \HH_j \otimes \hol_B$.  Now, by our calculations, $V_j = \sH_j$ over $B^* = B  \setminus S$, $S$ being the set of critical values of $f$.
We saw that the norm of a local frame of $V_j$ has at most  logarithmic grow at the  points $p \in S$.
This shows that $V_j $ is a subsheaf of $\sH_j$: by semipositivity we conclude that we have equality $V_j = \sH_j$.

\qed

\section{Counterexamples to Fujita's question}
In this section we shall provide two examples of surfaces fibred over a curve, with fibres curves with a symmetry of $ G : = \ZZ/7$ as in the preceding
section.

 We consider again the equation 
$$ z_1^7 = y_1 y_0 (y_1 - y_0 )  (y_1 - x y_0 )^4 , \ \  x \in \CC \setminus \{0,1  \}$$
but we homogenize it to obtain the equation

$$ z_1^7 = y_1 y_0 (y_1 - y_0 )  (x_0 y_1 - x_1 y_0 )^4  x_0^3. \  $$

The above equation describes a singular surface $\Sigma'$ which is a cyclic covering of 
 $\PP^1 \times \PP^1$ with group $ G : = \ZZ / 7$; $\Sigma'$ 
is contained inside the line bundle $\LL_1$  
over $\PP^1 \times \PP^1$  whose sheaf of holomorphic sections $\sL_1$ equals $ \hol_{\PP^1 \times \PP^1}(1,1)$.
One may observe that the second projection shows that the surface $\Sigma'$ is a ruled surface.

Since the branch divisor is a not a normal crossing divisor, we blow up the point $ x_0 = y_0 = 0$,
obtaining a del Pezzo surface which we denote by $Z$, while we denote by $\Sigma$ the normalization of
the induced $G$-Galois cover of $Z$. 

\begin{rem}
The singularities of the normal surface $\Sigma$ are of three analytical types, which we describe by their analytical equation 

\begin{enumerate}
\item
$ z^7 = x^4 y$: one for each singular fibre
\item
$ z^7 = x^3 y$: three on the fibre at infinity 
\item
$ z^7 = x y$: one on the fibre at infinity.

\end{enumerate}
\end{rem}

Finally, we let $Y$ be a minimal resolution of singularities of $\Sigma$.
Therefore $Y$ admits a fibration $\varphi : Y \ra \PP^1$ with fibres curves of genus $6$.

We let $X$ be the minimal resolution of the fibre product of $\varphi : Y \ra \PP^1$ with $ \psi : B \ra \PP^1$,
where $\psi$ is the $G$-Galois cover branched on $ \infty = \{x_0= 0 \}, 0 = \{ x_1 = 0\} , 1 = \{ x_1 = x_0\}$,
and with local characters $(1,1,-2)$. In particular $B$ has genus 3 by Hurwitz' formula ( $ 2g-2 =  7 \cdot ( -2 + 3 ( 1- \frac{1}{7}) \Rightarrow g = 3$).

Observe in fact that the singular fibres of $\varphi$ are exactly those lying above those three points.
Then there is a fibration  $ f : X \ra B$, with only three singular fibres.

We shall prove in a later subsection the following 

 \begin{theo}\label{firstsurface}
The above   surface $X$ is a surface of general type endowed with  a fibration $ f : X \ra B$ onto a curve $B$ of genus $ 3$, and with fibres of genus $6$,  such that $V : = f_* \om_{X|B}$ 
splits as a direct sum $ V = A  \oplus Q_1 \oplus Q_2$, where $A$ is an ample   rank-2  vector bundle, and the unitary  flat rank-2 summands $Q_1, Q_2$
have infinite monodromy group (i.e., the image of $\rho_j$ is infinite).
 \end{theo}

Consider now  the equation 
$$ z_1^7 = y_1 y_0^4 (y_1 - y_0 )  (y_1 - x y_0 ), \ \  x \in \CC \setminus \{0,1  \}$$
which gives another family of curves. It is similar to the previous family, 
except that we get here $V_1 $ generated by
$$ \eta: = y^{- \frac{6}{7}}  (y-1) ^{- \frac{6}{7}} (y-x) ^{- \frac{6}{7}} dy , {\rm and \ by \ } y \cdot \eta. $$ 

We shall see in the next section how, varying $x$, we obtain a rank-2 local system over $ \PP^1 \setminus \{ 0, 1, \infty\}$,
which is equivalent, in view of the Riemann-Hilbert correspondence,
to a second order differential equation with regular singular points.
Indeed,we shall see that we have in fact a Gauss hypergeometric equation.

But now we homogenize  the equation  to obtain

$$ z_1^7 = y_1 y_0^4 (y_1 - y_0 )  (x_0 y_1 - x_1 y_0 )  x_0^6. \  $$

This is a $G$-covering of $\PP^1 \times \PP^1$, and we
obtain another $G$-covering of $\PP^1 \times \PP^1$ by  taking its birational pull-back $T$

$$ z_1^7 = y_1 y_0^4 (y_1 - y_0 )  (P_7 (x)  y_1 - G_7 (x)  y_0 )  x_0^6,   $$
where $P_7, G_7$ are generic degree $7$ homogeneous polynomials.

We denote by $X(T)$ the minimal resolution of the singularities of $T$.

\begin{rem}
The singularities of the normal surface $T $ are of two analytical types, which we describe by their analytical equation 

\begin{enumerate}
\item
$ z^7 = y^4 u$: one for each point $ y_0 = P_7 (x) = 0$ (one point for seven fibres $ P_7 (x) = 0$)
\item
an $A_6$-singularity $ z^7 = u y$ for each point $ y_1 = G_7 (x) = 0$ (one point for seven fibres $ G_7 (x) = 0$),
and  for each point $ y_1- y_0  = ( P_7 + G_7) (x) = 0$ (one point for seven fibres $ ( P_7 + G_7) (x) = 0$).

\end{enumerate}
(1.1) A singularity $ z^7 = y^4 u$ is a quotient singularity of type $\frac{1}{7} (1,3)$: since, if we set $ x = w^7, y = v^7$,
then $ z : = v^4 w$ is invariant  for $ v \mapsto \z v,  w \mapsto \z^3 w$.
The minimal resolution of singularities is given by a Hirzebruch -Jung string of $\PP^1$'s $E_3, E_2, E'_2$ with respective self-intersections $ -3, -2 , -2$
(indeed $ \frac{7}{3} =  3 - \frac{1}{ 2 -  \frac{1}{ 2}}$).

(1.2) In our case, the fibre $\PP^1$ intersects transversally the two curves locally given by the equation $u=y= 0$.
Some calculations with the resolution of these quotient singularities (see \cite{bpv}, page 80) show that the
fibre of the minimal resolution $X(T)$ of $T$ consists of a smooth curve of genus three tangent to $E_2$ at the intersection point of $E_3$ and $E_2$.
We need two blow ups of this point to obtain that the fibre is a normal crossing divisor. Then the multiplicities
of the exceptional divisors are respectively $7,4,2,1,1$: we conclude then that in order to obtain the semistable reduction
we must take a covering of the base which is ramified at the point $P$ corresponding to the singular fibre of order divisible by 28.

(2.1) The $A_6$ singularity is resolved by a chain of $\PP^1$'s $$E_1, E_2, E_3, E_3', E'_2, E'_1$$ with self-intersection equal to $-2$.
The fibre of the minimal resolution  $X(T)$ of $T$ consists of a smooth curve of genus three intersecting $E_3$ and $E'_3$ transversally at the point
$E_3 \cap E'_3$, and the sum $(E_1+ E'_1) + 2 (E_2 +  E'_2) + 3  (E_3 +  E'_3)$. We need just to blow up
the point
$E_3 \cap E'_3$ to obtain a normal crossing divisor. Since the multiplicities of the seven exceptional divisors in the new chain are
$1,2,3,7,3,2,1,  $we conclude then that in order to obtain the semistable reduction
we must take a covering of the base which is ramified at the point $P$ corresponding to the singular fibre of order divisible by 42.

\end{rem}

\subsection{Associated local systems on $\PP^1$}
Let $P:=\PP^1,\, S:=\{s_1,\ldots,s_{r}\}$ and
 let $\LL$ be a rank-one 
local system 
on $P\setminus S$ corresponding to a homomorphism
$$ \rho : \pi_1(P)=\langle \gamma_{s_1}, \ldots,\gamma_{s_r} |  \gamma_{s_1}\cdot  \ldots \cdot \gamma_{s_r} = 1\rangle\To \CC^*,
 \gamma_s\Mapsto \alpha_s,$$
where $\gamma_s$ denotes a simple loop around $s\in S.$ We shall 
always
assume that, for all $s\in S,$ 
the monodromy generators $\alpha_s$
are roots of unity different from~$1.$ 

Observe that $\rho$ determines a Galois covering $\phi : C \ra P$ with Galois group $ G : = Im (\rho)$.
We have that $G = \mu_n : = \{ \zeta | \zeta ^n = 1\}.$
Hence we may write $\alpha_s=e^{2\pi i \frac{m_s}{n}},$
where $ 0 \leq m_s < n$, and we set also $\nu_s := \frac{m_s}{n}$.

The equation of $C$ is therefore given by
$$  z_1^n = \Pi_j (y_1 - s_j y_0)^{m_j}.$$ 

We have an eigenspace splitting for the direct image of the sheaf of holomorphic 1-forms:
$$ \phi_* (\Omega^1_C) =  \bigoplus_{h=0}^{n-1} \phi_* (\Omega^1_C)_h =  \bigoplus_{h=0}^{n-1} ( \Omega^1_P \otimes \sL_{-h}). $$ 

We want to relate the above summands to local systems on $P$. To this purpose, observe
that any character $$\chi_h
:\mu_n \To \CC^*, \zeta  \Mapsto \zeta^h$$ defines  a rank-one local system 
$\LL_j$ on $P\setminus S,$ associated to the homomorphism $\rho^h$.

Let $\LL$ be any of the $\LL_h$: then we have a Hodge decomposition
$$H^1(P\setminus S,\LL)=H^{(1,0)}(P\setminus S,\LL)\oplus H^{(0,1)}(P\setminus
S, \LL)$$ where $H^{(1,0)}(P\setminus S,\LL)$ is the space of  differentials 
of the first kind in $ H^0(P, j_*^{\rm mer}\Omega^1(\LL))$
and $ H^{(0,1)}(P\setminus
S, \LL)$ is  the complex conjugate (e.g., $\LL_{-j}$ is the complex conjugate of $\LL_j$)of the corresponding group for the 
dual (i.e., conjugate) local system
(cf.~\cite{DeligneMostow}, Page  19 and 
Prop. 2.20). 

Again by \cite{DeligneMostow}, Prop. 2.20,
 the Hermitian form on $H^1(P\setminus S,\LL)$ 
given by loc.cit., 2.18, (under the 
identification of $H^{1}_c(P\setminus
S,
\LL)$ with $H^1(P\setminus S,\LL)$) is positive definite on $H^{(1,0)}(P\setminus
S,\LL)$ and negative definite on $H^{(0,1)}(P\setminus
S,\LL).$ 

The relation between the two points of view is simply given by the equalities (compare \cite{DeligneMostow}):
$$ H^0  ( \Omega^1_C)_h =  H^{(1,0)}(P\setminus S,\LL_{-h}),$$

$$ H^1  ( C, \hol_C)_h = H^1  ( P, \sL_h^{-1}) = H^{(0,1)}(P\setminus S,\LL_{-h}).$$

Let $\nu_s \in \QQ$ as above be the unique rational number  between 
$0$ and $1,$ 
satisfying $\alpha_s=e^{2\pi i \nu_s}.$ By 
\cite{DeligneMostow}, Equation 2.20.1,
\begin{equation}\label{eq2}
 \dim H^{(1,0)}(P\setminus S,\LL)=-1+\sum_{s\in S}\mu_s.\end{equation}
 
In the case that the cardinality of $S$ is equal 
to $4$ 
we have $\dim H^1(P\setminus S,\LL)=2$ by \cite{DeligneMostow}, Prop. 2.3.1,  and 
$$ \dim(H^{(1,0)}(P\setminus S,\LL))=0,1,2,$$
corresponding to the cases 
$$ \sum_{s\in S}\mu_s=1,2,3,$$ by Formula~\eqref{eq2}. \\

Consider now  the  above family of projective 
curves $f: Y \to P$ birationally defined
by 
$$ z_1^7=y(y-1)(x-y)^4,$$
where $x$ denotes the affine coordinate of $P$ and let $C^o$ be its
restriction 
to $P\setminus \{0,1,\infty\}.$

Any character $$\chi_j
:\mu_7\To \CC^*, \zeta \Mapsto \zeta^j$$ defines  a rank-one local system 
$\LL_j$ on $P\setminus S,\, S=\{0,1,x,\infty\}$ with multi-valued local sections of the form
$ y^{j/7}(y-1)^{j/7}(y-x)^{4j/7}$ (\cite{DeligneMostow},  2.11),
having monodromy 
generators $\alpha_s=e^{2\pi i \frac{j}{7}}$ for $s=0,1,$
$\alpha_x=e^{2\pi i \frac{4j}{7}}$  and 
$\alpha_\infty=e^{-2\pi i \frac{(j+j+4j)}{7}}$.

It also gives rise to  the $\chi_j$-equivariant rank-two vector space 
$$ H^1(C_x,\CC)^{\chi_j}\simeq H^1(P\setminus S,\LL_{-j}).$$

The de Rham version of the cohomology group
$H^1(P\setminus S,\LL_6)$ is the fibre of a rank-2 vector bundle $E$ on 
$P\setminus \{0,1,\infty\}$  
with flat 
connection $\nabla:E\to \Omega^1_{P\setminus S}\otimes E$
whose local holomorphic solutions are integrals of the form 
$$g(x)=\int_1^{\infty} y^{-\frac{6}{7}}(y-1)^{-\frac{6}{7}}(x-y)^{-\frac{3}{7}}dy$$
(or similar integrals over Pochhammer double loops).

\subsection{Monodromy of some character sheaves}

By \cite{Kohno99}, Page 169, the above function $g(x)$ coincides 
up to a constant factor with the Gau{\ss} hypergeometric function
$$F(\alpha= 8/7,\beta= 3/7,\gamma= 9/7;x).$$ 
This implies that the rank two connection
$\nabla$ is equivalent to (the connection on $P\setminus \{0,1,\infty\}$  associated to) the 
Gau{\ss} hypergeometric differential
equation 
$$ t(t-1)f''+((\alpha+\beta+1)t-\gamma)f'+\alpha\beta f=0$$
(cf. \cite{Kohno99}, Page 163). 
The latter equation is non-resonant (i.e., the difference of two numbers of  $\{\alpha,\beta,\gamma\}$ 
does not lie in $\ZZ$), implying that the differential equation and hence its monodromy is irreducible.
It has the Riemann scheme
$$ \left\{ \begin{array}{ccc}
0&1&\infty\\
0&0&\alpha\\
1-\gamma& \gamma-\alpha-\beta&\beta
\end{array}\right\},$$
cf. \cite{Kohno99}, Page 164. Since the Riemann scheme 
describes the  exponents of a basis of solutions
of an ordinary differential equation in their respective 
Puiseux expansions, 
this  implies that the local monodromy of $\nabla$
at $0,1$ is a homology of order $7$ and hence is of order 
$7$ in the associated projective linear group. Recall the 
Schwarz' list of the Gau{\ss} hypergeometric differential equations
 with finite projective monodromy groups
\cite{Schwarz}. As no two consecutive projective 
local monodromies of order $7$ occur amongst the irreducible cases listed there, 
 we conclude that the monodromy of $\nabla$ 
is infinite.\\

\subsection{Proof of theorem \ref{firstsurface}} 

We considered a ramified covering $\psi :B\to  P$ which 
is locally at each branch point $0,1,\infty$ 
of type $x\mapsto x^7.$
Then we got $ f : X \ra B$ as the minimal resolution of 
the fibre product $B\times_P Y \to B.$

The fibres of $f$ are smooth curves of genus 6 and $X$ has an action of $ G \cong \mu_7$
which is of type $(4,1,1,1)$ on all smooth fibres.

There are only three singular fibres, but around them the monodromy
of the rank-2 local systems $\HH_1^*, \HH_2^*$ is trivial, because we saw in the previous section
that the local monodromy is of order 7.

Hence these extend to rank-2 local systems $\HH_1, \HH_2$ over $B$.

The same argument given in proposition \ref{7} shows then that $V_j = \sH_j$
for $j=1,2$. We have then $ V = U \oplus Q_1 \oplus Q_2$, where we set  $Q_j : = V_j = \sH_j$
for $j=1,2$, and $ U : = V_3 \oplus V_4$.

Assume that $U$ is not ample, and that it contains a unitary  flat summand $Q'$. 

Without loss of generality, we may
assume that $ Q' | B^*\subset \HH_3^*$.

Since $\HH_3^* = V_3 \oplus \overline{V_4}| B^*$ we see that $Q'$ has rank $1$.

By the cited result by Deligne ((I) in corollary \ref{cor}) $Q'$ would be a torsion line bundle, hence also $V_3$ and $V_4$
and the monodromy  of  $\HH_3^*$ (respectively $\HH_4^*$) would be finite.

However, the integrals associated to the factor $\HH_3^*$ (respectively $\HH_4^*$) also satisfy a Gauss
differential equation with infinite  monodromy (again by Schwarz' list, since the local monodromies  at $0,1,\infty \in P$ are of order $\geq 7$):
this gives a contradiction.

\begin{rem}The same considerations apply to the second family of curves that we introduced,
and also  to other families of 
curves.
\end{rem}

\bigskip

\noindent
{\bf Acknowledgements.}
The interest of the first author in the Fujita question was aroused by the work of Miguel Angel Barja (\cite{Barja}).

We thank also  Stefan Reiter for valuable conversations on hypergeometric differential equations, and  Thomas Peternell and Yujiro Kawamata 
for a useful question which helped us to improve the exposition.


\end{document}